\documentclass{amsart}
\usepackage{amssymb}
\usepackage{amsfonts, a4wide}
\usepackage{amsmath, amsthm}

\usepackage[all]{xy}

\newtheorem{Theorem}{Theorem}[section]
\newtheorem{Lemma}[Theorem]{\sc Lemma}

\newtheorem{Question}[Theorem]{\sc Question}

\def\Q{{\mathbb{Q}}}
\def\C{{\mathbb{C}}}
\def\R{{\mathbb{R}}}
\def\Z{{\mathbb{Z}}}

\def\pf{{\noindent{\sc Proof: }}}
\def\F{\mathbb{F}}
\def\urat{{\widetilde{Rat}}}
\def\uc{{\widetilde{C}}}

\newcommand\ra{\rightarrow}
\newcommand\la{\leftarrow}
\newcommand{\bt}{\begin{Theorem}}
\newcommand{\et}{\end{Theorem}}
\newcommand{\bl}{\begin{Lemma}}
\newcommand{\el}{\end{Lemma}}
\newcommand{\bd}{\begin{displaymath}}
\newcommand{\ed}{\end{displaymath}}
\newcommand{\bq}{\begin{Question}}
\newcommand{\eq}{\end{Question}}
\newcommand{\ba}{\begin{align*}}
\newcommand{\ea}{\end{align*}}

\begin{document}

\title{The Cohomology Ring of the Space of Rational Functions}
\author[Dinesh Deshpande]{Dinesh Deshpande}
\address{Department of Pure Mathematics and Mathematical Statistics, \\University of Cambridge, \\Cambridge CB3 0WB, UK}
\email{D.Deshpande@dpmms.cam.ac.uk}

\maketitle

\begin{abstract}
 Let $Rat_k$ be the space of based holomorphic maps from $S^2$ to
 itself of degree $k$. Let $\beta_k$ denote the Artin's braid group on
 $k$ strings and let 
$B\beta_k$ be the classifying space of $\beta_k$. Let $C_k$
denote the space of configurations of length less than or equal to $k$
of distinct points in $\R^2$ with labels in $S^1$.  The three spaces
$Rat_k$, $B\beta_{2k}$, $C_k$ 
 are all stably homotopy equivalent to each other. For an odd prime $p$, the
$\F_p$-cohomology ring of the three spaces are isomorphic to each
other. The $\F_2$-cohomology ring of $B\beta_{2k}$ is 
isomorphic to that of $C_k$.  
 We show that for all values of $k$ except 1
and 3, the $\F_2$-cohomology ring of $Rat_k$ is not isomorphic to that of
$B\beta_{2k}$ or $C_k$. This in particular implies that the
$H\F_2$-localization of  $Rat_k$ is not
homotopy equivalent to $H\F_2$-localization of $B\beta_{2k}$ or
$C_k$. We also show that for 
$k \geq 1 $, $B\beta_{2k}$ and $B\beta_{2k+1}$ have homotopy equivalent $H\F_2$-localizations.
\end{abstract}

\section{Introduction}

Let $\beta_k$ denote the Artin's braid group on $k$ strings. Let
$B\beta_k$ be the classifying space of $\beta_k$. Let $C_k(\R^2,S^1)$
denote the space of configurations of length less than or equal to $k$
of distict points in $\R^2$ with labels in $S^1$, with some
identifications. We use just $C_k$ to denote the space
$C_k(\R^2,S^1)$.  Let $Rat_k$ be the
space of based holomorphic maps from $S^2$ to itself of degree $k$.   

 \cite{seg79},\cite{cohcoh91,cohcoh93} shows that theses
three spaces $Rat_k$, $B\beta_{2k}$, $C_k$ are all stably homotopy
equivalent. In fact, \cite{cohdav88, sna74} shows that these spaces split
stably as a wedge sum $\vee_{j \leq k} D_j(S^1)$, where $D_j =
C_j/C_{j-1}$ is a space related to the Brown-Gitler spectra. The three
spaces are closely related to $\Omega^2S^2$. We explain some facts about these spaces in the next section. 
 
Totaro \cite{tot90} has shown that the three spaces have isomorphic
$\F_p$-cohomologies for an odd prime $p$. He has also shown that the
$\F_2$-cohomology ring of $B\beta_{2k}$ is isomorphic to that of $C_k$ and
if $k+1$ is not a power of 2, then the $\F_2$-cohomology ring of $Rat_k$ is not isomorphic to that of $B\beta_{2k}$ or $C_k$. 
This paper extends the result
to all values of $k$ except when $k = 1$ or $k = 3$ [Theorem
\ref{main}]. This in particular implies that $Rat_k$ is not homotopy equivalent to $C_k$ if $k$ is not equal to 1 or 3. Bousfield has
defined the localization of spaces with respect to homology in
\cite{bou75}. Two spaces $X$ and $Y$ have homotopy equivalent
$HR$-localizations if and only if there are maps $X \ra X_1 \la Y$  such that
each map induces an isomorphism on homology groups with coefficients
in the ring $R = \Z, \F_p$ or $\Z[q^{-1}]$. Our result implies that,
$H\F_2$-localizations of $Rat_k$ and $B\beta_{2k}$ are not homotopy equivalent. 
 We also show that
for $k \geq 1$ 
$B\beta_{2k}$ and $B\beta_{2k+1}$ have isomorphic
$\F_2$-cohomologies and $H\F_2$-localizations of $B\beta_{2k}$ and
$B\beta_{2k+1}$ are homotopy equivalent [Lemma \ref{braid}]. 

For $k=1$, the three
spaces $Rat_1$, $B\beta_2$ and $C_1$ are all homotopy equivalent to
$S^1$. For $k = 3$, it turns out that the corresponding three spaces have
isomorphic cohomology rings with coefficients in $\F_p$ for any prime
$p$. Moreover, the actions of the dual of the Steenrod algebra on the $\F_2$-homologies of $Rat_3$, $B\beta_6$ and $C_3$ are also isomorphic.\newline

{\it Acknowledgements}: The author thanks his PhD supervisor Burt Totaro for introducing to this subject and for numerous interesting discussions.

\section{$B\beta_k$, $C_k$, $Rat_k$}
In 
this  section, we describe the three
spaces $Rat_k$, $B\beta_k$ and $C_k$, their respective integral cohomologies and
their relation with the space $\Omega^2S^2$. We also describe the 
coalgebra structure of 
their respective
$\F_2$-homologies. For this chapter, the default ring of coefficients is
$\F_2$.

  Let $\Omega^2S^2$ be the
double loop space of $S^2$, i.e the space of maps from $S^2$ to itself.
\begin{align*}
 \pi_0(\Omega^2S^2) \cong \pi_2(S^2) \cong \Z, 
\end{align*}

where the degree map
induces an isomorphism
\begin{align*}
\xymatrix{
  \pi_0(\Omega^2S^2) \ar[r]^{\cong} & \Z.}
\end{align*}

\subsection{$B\beta_k$}

 The braid space $B\beta_k$ is the classifying space of the braid group
on $k$-strings, $\beta_k$. Let $F(\R^2, k)$ denote the
configuration space of $k$-points in $\R^2$, i.e. \bd F(\R^2, k) = \{
(x_1,...,x_k)| x_i \in \R^2,i \neq j \Rightarrow  x_i \neq x_j
\}.\ed The symmetric group on $k$ elements $\Sigma_k$ acts freely on
$F(\R^2, k)$. We can take $ F(\R^2, k)/ \Sigma_k$ as a model of the
classifying space of $\beta_k$. Thus the space
$B\beta_k$ is the
space of unordered $k$-tuples of points in $\R^2$. The space
$B\beta_k$ can also be described as the space of degree $k$ complex
polynomials without multiple roots and with the leading coefficient
equal to the unity.  

The rational cohomology of the braid groups is as follows
(\cite{ver99}, Theorem 8.1-2).

\bl For $k \geq 2$, the rational cohomology groups of $B\beta_k$ are
trivial except for
\begin{align*}
  H^0(B\beta_k;\Q) \cong& \Q,\\
H^1(B\beta_k;\Q) \cong& \Q.
\end{align*}

And for $k \geq 1$,
\begin{align*}
  H^i(B\beta_{2k+1};\Z) \cong H^i(B\beta_{2k};\Z).
\end{align*}

\el

As the spaces $B\beta_{2k}$, $Rat_k$ and $C_k$ are stably homotopy equivalent to each other,
\begin{align*}
 H^i(B\beta_{2k};\Z) \cong H^i(Rat_k;\Z) \cong H^i(C_k;\Z).  
\end{align*}

F. Cohen has calculated the $\F_p$-homology of
$\coprod B\beta_k$ and 
$\Omega^2S^2$ in \cite{cohlad76}.  The 
spaces $\Omega^2S^2$ and $B\beta_k$ are $\mathcal{C}_2$-spaces,
i.e. $\mathcal{C}_2$, the `little 2-cubes operad' acts on them. Hence there  
is the Araki-Kudo operation on the $\F_2$-homologies of
$\coprod_k B\beta_k$ and $\Omega^2S^2$,  
$ Q : H_q \ra H_{2q+1}  $ and the Pontrjagin product   which
makes their homologies commutative rings. 
Let $\Omega_k^2 S^2$ denote the $k^{th}$
component of $\Omega^2S^2$ corresponding to the degree $k$ maps $f: S^2
\ra S^2$. Then $Q$ maps $H_i(\Omega_k^2S^2)$ to $H_{2i+1}(\Omega^2_{2k}S^2)$.

 There is a natural map $\phi: B\beta_k \ra \Omega^2_kS^2$. This map can be
 described as follows. Replace the 
 $k$-tuple of distinct points in $\R^2$ by $k$ disjoint unit circles in
 $\R^2$. Then define a map from $\R^2 \cup \infty$ to itself 
 by sending everything except interiors of unit
 circles to the point at infinity and by sending the interior of each
 unit circle onto the whole of $\R^2$ homeomorphically. Identifying $\R^2
 \cup \infty$ with $S^2$ by the stereographic projection gives a degree
 $k$ map from $S^2$
 to itself. This is precisely the natural map $\phi$ from
 $B\beta_k$ to $\Omega^2_kS^2$. An algebraic construction of a map
 $B\beta_k \ra \Omega^2S^2$ is given in section 1, \cite{seg79}.

Note that
\begin{align*}
  \pi_1(B\beta_k)& \cong \beta_k,\\
\pi_1(\Omega^2S^2)& \cong \Z,\\
\pi_n(B\beta_k)& \cong \{0\}, \forall k > 1.
\end{align*}

Hence the map $\phi$ can not be a homotopy equivalence in any range of
dimensions. But it turns out that the map induces an isomorphism of
homologies up to dimension $\lfloor k/2 \rfloor := $\emph{the greatest integer
  smaller than or equal to }$k/2$.

 The map $\phi$ 
induces a map $\Phi: \coprod_{k \geq 0} B\beta_k  \ra \Omega^2S^2$. Let $g$ be the
generator of $H_0(B\beta_1)$. By using the map $\Phi$, let
$g$ also denote the generator of $H_0(\Omega^2_1S^2)$. Then the homology of
these two spaces is build-up by the `Araki-Kudo' operation $Q$ and its
iterations $Q^i(x) = Q(Q^{i-1}(x))$. To be
precise, there are algebra isomorphisms (appendix III,
\cite{cohlad76}) \bd H_*(\coprod B\beta_k) \cong 
\F_2[g, Qg, Q^2g, \ldots], \ed \bd H_*(\Omega^2S^2) \cong \F_2[g,
g^{-1}] \otimes \F_2[Qg, Q^2g, \ldots].  \ed

Note that the dimension in homology of $Q^ig$ is  $2^i-1$ and is
contained in the
$2^i$th component of $\Omega^2S^2$. Define the weight of a homology
class to be the component in which that class lives. Hence,
$H_*(B\beta_k)$ is the span of monomials in $g, Qg, Q^2g, ...$ of
weight $k$, where $Q^ig$ has the weight $2^i$ and the dimension $2^i-1$. Hence
note that for any $k$, the top dimensional homology of $B\beta_k$ is
generated by a single element. If the binary expansion of $k$ is 
$k = \sum_{j \in J} 2^j$, then this top dimension is $H_{k-|J|}$. Also
notice that $Q(x^2) = x^2Qx + Qx \cdot x^2 = 0$ as the homology
coefficients are in $\F_2$.

Further,  this operation $Q$ is linear and that the Cartan formula
holds (lemma 5.2, IX, \cite{cohlad76}) \bd Q(xy) = x^2Qy + Qx \cdot y^2. \ed

The coproduct structure on the homology, i.e. the cup product strucutre on the
cohomology of $B\beta_k$ is as given below. It turns out that
$H_*(B\beta_k)$ is a primitively generated Hopf algebra. i.e., let
\bd \psi: H_* \ra H_* \otimes H_* \ed denote the coproduct on the
homology. Then $\psi(g) = g \otimes g$ and $Q^ig$ for $i \geq 1$ is
primitive in its component, \bd \psi(Q^ig) = g^{2^i} \otimes Q^ig +
Q^ig \otimes g^{2^i}. \ed 

$\psi$ being a coproduct map satisfies that \bd \psi(xy) = \psi(x)\psi(y).  \ed
The expressions for $Q(g^{-1})$ and $Q(g^{-1}Qg)$ in
$H_*(\Omega^2S^2)$ in terms of $g$ and
$Q^ig$ can be obtained using the Cartan formula. They are,

\begin{align*} Q(g^{-1}) &= g^{-4}Q(g) \\ Q(g^{-1}Qg) &= g^{-2}Q^2g +
  g^{-4}(Qg)^3  \end{align*}

$H^*(B\beta_{2k})$ is isomorphic to $H^*(B\beta_{2k+1})$ \cite{arn68,fuk70}. We will show that $B\beta_{2k}$ and $B\beta_{2k+1}$ have
homotopy equivalent $H\F_2$-localizations.

\bl\label{braid} The cohomology ring  $H^*(B\beta_{2k})$ is isomorphic to
$H^*(B\beta_{2k+1})$. In fact, $B\beta_{2k}$ and $B\beta_{2k+1}$ have
homotopy equivalent $H\F_2$-localizations.
\el
\pf Let $x \in H_i(B\beta_{2k})$. Hence $x$ is an element of weight
$2k$ and dimension $i$ in $\Z_2[g,Qg,Q^2g,\cdots]$. Hence $gx$ is an
element of weight $2k+1$ and dimension $i$. Hence $gx \in
H_i(B\beta_{2k+1})$. Also, let  $y$ be a monomial in
$\Z_2[g,Qg,Q^2g,\cdots]$ of the
dimension $i$ and the weight $2k+1$, i.e. $y \in H_i(B\beta_{2k+1})
$. As each of the $Q^ig$ has even 
weight, $y$ is divisible by $g$, and $y/g = x \in H_i(B\beta_{2k})$. Also
\begin{align*}
  \psi(y) =& \psi(g)\psi(x)\\
=& (g \otimes g)\psi(x).
\end{align*}
 
Hence multiplication by $g$ induces an isomorphism of coalgebras
\begin{align*}
  \xymatrix{
H_*(B\beta_{2k}) \ar[r]^{\cdot g} & H_*(B\beta_{2k+1}).}
\end{align*}

Hence,
\begin{align*}
  H^*(B\beta_{2k}) \cong H^*(B\beta_{2k+1}).
\end{align*}

Furthermore, let $i_{2k}: B\beta_{2k} \ra B\beta_{2k+1}$ be the
inclusion map given by adding a point away from a given $2k$-tuple to
get a $2k+1$-tuple. Then note that $i_{2k_*}$ is precisely
multiplication by $g$. Hence $i_{2k_*}$ induces isomorphism on
$\F_2$-homologies. Hence $B\beta_{2k}$ and $B\beta_{2k+1}$ have
homotopy equivalent $H\F_2$-localizations.

\qed
 
The action of the dual of the Steenrod algebra on $H_*(\coprod B\beta_k)$ is given in
the appendix of \cite{coh78}. Let $Sq_j^*: H_n(-) \ra H_{n-j}(-)$ be the dual of the $j$th Steenrod operation $Sq^j$. 
Then 
\begin{align*} Sq_j^*(Q^ig) =& 0 \text{ if } j \geq 2.\\
Sq_1^*(Q^ig) =& (Q^{i-1}g)^2 \text{ if } i \geq 2\\
Sq_1^*(Qg) =& 0. \end{align*}

\subsection{$C_k$} 

The following description of the configuration spaces is taken from
\cite{tot90}, 
i.e. in turn from \cite{cohshi91}. 

Let $C(\R^2,Y)$ denote the space of all configurations of distinct points 
in $\R^2$ with labels in $Y$. It is defined by \bd C(\R^2,Y) =
(\bigcup_{j=1}^{\infty} F(\R^2,j) \times_{\Sigma_j} Y^j)/\sim \ed
and if $* \in Y$ is a fixed basepoint then
the equivalence relation $\sim$  is given by \bd (x_1,\ldots, x_j)
\times_{\Sigma_j} (t_1,\ldots , t_{j-1},*) \sim (x_1,\ldots,x_{j-1})
\times_{\Sigma_{j-1}}(t_1,\ldots,t_{j-1}).  \ed 

Let $C_k(\R^2,Y)$ denote the subspace of all configurations of length less
than or equal to $k$. i.e. \bd C_k(\R^2,Y) = (\bigcup_{j=1}^k F(\R^2,j)
\times_{\Sigma_j} Y^j)/\sim . \ed

We denote by  $C_k$ the space $C_k(\R^2, S^1).$

There is a  relation between configuration spaces and iterated loop
spaces(May-Milgram and Segal ). 
If $Y$ is a connected CW-complex then $C(\R^2,Y)$ 
is homotopy equivalent to the based loop
space $\Omega^2\Sigma^2Y$ which is defined by \bd \Omega^2\Sigma^2Y =
\{f: S^2 \ra \Sigma^2Y | f(\infty) = *  \}  \ed 
Hence $C_k$ can be considered as a finite dimensional approximation to
$\Omega^2S^3$. \bd \pi_1(C_k) \cong \Z.  \ed The Hopf map $S^3 \ra
S^2$ induces a map of 2-fold loop 
spaces, from $\Omega^2S^3$ to $\Omega^2S^2$. The long exact
sequence of the homotopy groups of the fibration $S^1 \ra S^3 \ra
S^2$ implies that $\Omega^2S^3 \ra  \Omega^2S^2$ gives the homotopy
equivalence from $\Omega^2S^3$ to $\Omega^2_0S^2$. This helps in
obtaining the following result (theorem 3.1, III, \cite{cohlad76}),

\bd H_*(\bigcup_{k \geq 0}C_k) \cong H_*(\Omega_0^2S^2) \cong
\F_2[g^{-2}Qg, Q(g^{-2}Qg), \ldots]. \ed

$H_*C_k$ is the span of monomials of weight less than or equal to $k$,
where the weight of $Q^i(g^{-2}Qg)$ is $2^i$ and it lives in the dimension
$2^{i+1} - 1$. Proposition 1 from \cite{tot90} shows that as coalgebras,
\bd H_*(B\beta_{2k}) \cong H_*(C_k).  \ed

 Havlicek \cite{hav95} has described the precise cohomology ring of
$C_k$ as the dual to this coalgebra.
  
\subsection{$Rat_k$}

The space $Rat_k(\C P^1)$ or $Rat_k$ is the space of based holomorphic
maps $S^2 \ra S^2$ of degree $k$. It can be described more precisely
as the space of rational functions from $\C \cup \infty$ to 
$\C \cup \infty$ which
 sends $\infty$ to 1, i.e. \bd Rat_k := \{
\dfrac{f(z)}{h(z)} = \dfrac{z^k + a_{k-1}z^{k-1} + \ldots + a_0}{z^k +
  b_{k-1}z^{k-1} + \ldots + b_0}| \text{\emph{f}(\emph{z}) and
  \emph{h}(\emph{z}) are coprime} \}  \ed   
 
$Rat_k$ is a nilpotent space up to dimension $k$ (corollary 6.3,
\cite{seg79}). i.e. the fundamental group of $Rat_k$ acts nilpotently
on homotopy groups $\pi_i(Rat_k)$ for $2 \leq i \leq k$. Consider the
map given by resultant of two polynomials  
\begin{align*}
  R : Rat_k \ra& \C^*\\
 (f/h) \mapsto& resultant(f,h).
\end{align*}

Then the map $R$ induces an isomorphism of fundamental groups
(proposition 6.4, \cite{seg79}) \bd \pi_1(Rat_k) \cong \Z.  \ed

There is a natural map $ Rat_k \ra
\Omega^2_kS^2$ which simply forgets that a map in $Rat_k$ is
holomorphic. This map is well described in (\cite{boyman88}, 
\cite{seg79}). This induces a map
\begin{align*}
  \chi: \coprod_{k \geq 0} Rat_k \ra \Omega^2S^2.
\end{align*}

The map $\chi$ preserves the action of $\mathcal{C}_2$ operad on the
spaces $\coprod_{k \geq 0} Rat_k$ and $\Omega^2S^2$.

The map $\chi$ induces a  
map on the homologies and in the proof of Theorem1 in \cite{tot90}, it is
shown that this induced map is an injection. The image of this map is a
polynomial ring generated by $g$ and $Q^i(g^{-1}Qg)$ for $i \geq
0$. To be precise, 
\bd H_*(\coprod_k Rat_k) = \F_2[g, g^{-1}Qg, Q(g^{-1}Qg), \ldots].
\ed  

As before, $g$ has weight 1 and dimension zero, and
$Q^{i-1}(g^{-1}Qg)$, $i \geq 1$ has weight $2^{i-1}$ and dimension
$2^i - 1$. $H_*(Rat_k)$ as a sub-coalgebra of $H_*(\coprod_k Rat_k)$ is
generated by the monomials of weight $k$. Note
again that the top dimensional homology of $Rat_k$ is generated by a
single element.

There is a one-to-one correspondence between the generators of
$H_*(\coprod_k Rat_k)$ and $H_*(\coprod_k B\beta_{k})$. $g$ of course 
corresponds to $g$ and $Q^{i-1}(g^{-1}Qg)$ corresponds to $Q^ig$. Note 
that in this correspondence, except for $g$, the weights of the generators of 
$H_*(\coprod_k Rat_k)$ are exactly the half of the weights of the generators
they correspond to in  $H_*(\coprod_k B\beta_k)$. 

The coproduct structure on $H_*(Rat_k)$  is as 
follows. $g^{-1}Qg$ is primitive in its component, but $Q^i(g^{-1}Qg)$  
is not primitive in its component for $i \geq 1$. 

\begin{align*} \psi(g^{-1}Qg) = g \otimes g^{-1}Qg + g^{-1}Qg \otimes g.
   \end{align*}
And for $i \geq 1$,
\begin{align*} \psi Q^i(g^{-1}Qg) = &\underbrace{g^{2^i} \otimes
    Q^i(g^{-1}Qg)}_{0,2^{i+1}-1} + \underbrace{Q^ig \otimes
    (g^{-1}Qg)^{2^i}}_{2^i-1,2^i} + \underbrace{(g^{-1}Qg)^{2^i}
    \otimes Q^ig}_{2^i,2^i-1} \\ &+ \underbrace{Q^i(g^{-1}Qg) \otimes
    g^{2^i}}_{2^{i+1}-1,0}.  \end{align*} 

Numbers appearing below a symbol indicate the dimension in homology of the
corresponding symbol. i.e. $(0,2^{i+1}-1)$ below $g^{2^{i}} \otimes
Q^i(g^{-1}Qg)$ indicate that $g^{2^i}$ is zero dimensional and
dimension of $Q^i(g^{-1}Qg)$ is $2^{i+1}-1$.

\section{$\F_2$-cohomology ring of $Rat_k$ is not isomorphic to that of
  $B\beta_{2k}$ or $C_k$}

This section proves that if $k$ is not equal to 1 or 3, then the
cohomology ring of $Rat_k$ with coefficients in $F_2$ is not
isomorphic to the cohomology ring of $B\beta_{2k}$ or $C_k$ with coefficients
in $F_2$.
 Totaro \cite{tot90} has shown this statement when $k+1$ is
not a power of 2. For $k=1$, all three spaces $Rat_1$, $B\beta_2$ and
$C_1$ are homotopy equivalent ot $S^1$. For $k=3$, the three spaces
$Rat_3$, $B\beta_6$ and $C_3$ have isomorphic $\F_p$-homology as  
coalgebras for any prime $p$. 
Following theorem shows the result for the remaining values of
$k$, that is when $k+1$ is a power of 2 and $k>3$.  This in particular
shows that if $k$ is not equal to 1 or 3, then there does not exists
any sequence of maps
\begin{align*}
  B\beta_{2k} \ra X_1 \la Rat_k
\end{align*}
each of which induces isomorphism on $\F_2$-homology. Hence in the context of
\cite{bou75} our result implies that the two spaces $Rat_k$ and
$B\beta_{2k}$ can not have homotopy equivalent $HF_2$-localizations. For
completeness, we
include Totaro's argument when $k+1$ is not a power of 2.

\begin{Theorem}\label{main} The $\F_2$-cohomology of 
 $Rat_k$ is not
  isomorphic to the $\F_2$-cohomology of $B\beta_{2k}$ or $C_k$ except
  when $k = 1$  or $3$. 
\end{Theorem} 
\pf

Firstly assume that $k+1$ is not a power of 2. Let
\begin{align*}
  k = \sum_{j \in J} 2^j.
\end{align*}

The top dimensional homology group of both $Rat_k$ and $B\beta_{2k}$ is 
1-dimensional. $H_{2k-|J|}$, the top dimensional homology of $Rat_k$ is spanned by
$x$ equal to $\prod_{j\in J} Q^j(g^{-1}Qg)$ and of $B\beta_{2k}$ by $y$ equal to
$\prod_{j \in J} Q^{j+1}g$.

Consider the set \bd S(x) = \{s \geq 0 | \psi(x)|_{H_s \otimes H_{d-s}} \neq 0
\}.\ed Similarly define $S(y)$. The aim is to show that $S(x)$ is not equal to
$S(y)$ which implies that the homologies of $Rat_k$ and $B\beta_{2k}$ are not
 isomorphic as coalgebras. 

Let $r$ be the smallest integer such that $r \in J$ but $r-1 \notin
J$. As $k+1$ is not a power of 2, such a $r$ exists. We observe that
$\psi(x)$ is non-zero in $H_{2^r-1} \otimes H_{dim(x)-(2^r-1)}$.
\begin{align*}
  \psi(x) = \prod_{j \in J} \psi(Q^j(g^{-1}Qg)).
\end{align*}

There is at least one term in
dimension $H_{2^r-1} \otimes H_{dim(x)-(2^r-1)}$ in the expansion of $ \psi(Q^j(g^{-1}Qg))$, which is
\begin{align*}
  Q^rg\prod_{j\in J, j \neq r}g^{2^j} \otimes (g^{-1}Qg)^{2^r}\prod_{j
    \in J, j \neq r} Q^j(g^{-1}Qg).
\end{align*}

As $r-1 \notin J$, there is no other term of this dimension in
$\psi(x)$. Hence $2^r-1 \in S(x)$.  

\bd \psi(y)  = \prod_{j \in J}   \Big( \underbrace{g^{2^{j+1}} \otimes
      Q^{j+1}g}_{0,2^{j+1}-1} + \underbrace{Q^{j+1}g \otimes
      g^{2^{j+1}}}_{2^{j+1}-1,0} \Big). \ed  

Note that $\psi(y)$ is zero in dimension $H_{2^r-1}
\otimes H_{dim(x) - (2^r-1)}$ as $r-1 \notin J$. Hence $2^{r}-1 \notin
S(y)$ and $S(x) \neq
S(y)$. This proves that whenever $k+1$ is not a power of 2, the
cohomology rings $H^*(Rat_k)$ and $H^*(B\beta_{2k})$ or $H^*(C_k)$ are
not isomorphic.

Now assume that $k+1$ is a power of 2, and assume that \bd k =
\sum_0^r 2^j = 2^{r+1} -1.\ed 

We continue to denote by $x$ the generator of the top dimensional homology 
group  of $Rat_k$, \bd x = \prod_0^{r} Q^j(q^{-1}Qg) \ed and by $y$,
the generator of the top dimensional homology group of $B\beta_{2k}$,
\bd y = \prod_0^{r} Q^{j+1}g.  \ed
For both spaces, the top dimension of homology is $d = 2^{r+2} - r
-3$. 

Let $S(x)$ and $S(y)$ be as before and we
will show that $S(x) \neq S(y)$.

\begin{align*}
\psi(y) &= \prod_0^{r} \psi(Q^{j+1}g)\\
 &= \prod_0^r \Big( \underbrace{g^{2^{j+1}} \otimes
      Q^{j+1}g}_{0,2^{j+1}-1} + \underbrace{Q^{j+1}g \otimes
      g^{2^{j+1}}}_{2^{j+1}-1,0} \Big) 
\end{align*}

Numbers $(0,2^{j+1}-1)$ appearing below $g^{2^{j+1}} \otimes
Q^{j+1}g$ indicate the dimension in homology of the corresponding element,
i.e. $g^{2^{j+1}}$ has dimension $0$ and $Q^{j+1}g$ has dimension
$2^{j+1}-1$. Hence dimensions which appear in $\psi(y)$
are precisely those which appear in the expression \bd \Big( ((1,0) +
(0,1)) \cdot ((0,3) + (3,0)) \cdot ((0,7)+(7,0)) \cdot \ldots
\Big) \ed 
From this expression, it is clear that $2 \notin S(y)$ and $5 \notin
S(y)$.

It turns out that although $2 \notin S(x)$, $5 \in S(x)$.

\bd \psi(x) = \psi(g^{-1}Qg)\prod_1^r \psi(Q^j(g^{-1}Qg)) \ed 

Using the expressions for $\psi(g^{-1}Qg)$ and $\psi(Q^j(g^{-1}Qg))$, we
get that the dimensions which appear in $\psi(x)$ are from the expression
\bd ((0,1)+(1,0)) \cdot \Big( ((0,3)+(1,2)+(2,1)+(3,0)) \cdot
((0,7)+(3,4)+(4,3)+(7,0)) \cdot \ldots   \Big)  \ed 

There are exactly two ways to obtain dimension $(2,d-2)$. Namely, $(1,0)
\otimes (1,2) \otimes (0,7) \otimes \ldots$ and $(0,1) \otimes (2,1)
\otimes (0,7) \otimes \ldots$. The expression corresponding to dimensions
$(1,0) \otimes (1,2)$ is $(g^{-1}Qg \otimes g) \cdot (Qg \otimes
(g^{-1}Qg)^2)$, which is $g^{-1}(Qg)^2 \otimes g^{-1}(Qg)^2$. By
symmetry, expression corresponding to $(0,1) \otimes (2,1)$ is also 
$g^{-1}(Qg)^2 \otimes g^{-1}(Qg)^2$. Hence, 
\begin{align*} \psi(x)|_{H_2 \otimes H_{d-2}} &= g^{-1}(Qg)^2 \otimes
g^{-1}(Qg)^2 \prod_2^r (g^{2^j} \otimes Q^j(g^{-1}Qg))\\ &+ g^{-1}(Qg)^2 \otimes
g^{-1}(Qg)^2 \prod_2^r (g^{2^j} \otimes Q^j(g^{-1}Qg)). \end{align*}
Because the coefficients of the homology  are  in $\F_2$, the two terms cancel
each  other. Hence, $2 \notin S(x)$.

There are exactly four ways to obtain dimension $(5,d-5)$. Namely
$(1,0) \otimes (1,2) \otimes (3,4) $, $(0,1) \otimes (2,1) \otimes
(3,4)$, $(1,0) \otimes (0,3) \otimes (4,3)$ and $(0,1) \otimes (1,2) \otimes
(4,3)$. From the paragraph above, first two of these, $(1,0) \otimes (1,2)
\otimes (3,4) $ and $(0,1) \otimes (2,1) \otimes (3,4)$ cancel with
each other. Hence,
\begin{align*}\psi(x)|_{H_5 \otimes H_{d-5}} &= (g^{-1}Qg \otimes
  g)(g^2 \otimes Q(g^{-1}Qg))((g^{-1}Qg)^4 \otimes Q^2(g)) \prod_3^r
  (g^{2^{j}} \otimes Q^j(g^{-1}Qg))\\ &+ (g \otimes g^{-1}Qg)(Qg \otimes
  (g^{-1}Qg)^2)((g^{-1}Qg)^4 \otimes Q^2g) \prod_3^r (g^{2^{j}} \otimes
  Q^j(g^{-1}Qg)) \end{align*}

\begin{align*} \Rightarrow \psi(x)|_{H_5 \otimes H_{d-5}} &=
  g^2(g^{1}Qg)^5 \otimes gQ^2gQ(g^{-1}Qg) \prod_3^r (g^{2^{j}} \otimes
  Q^j(g^{-1}Qg))  \\ &+ g^2(g^{1}Qg)^5 \otimes
  (g^{-1}Qg)^3Q^2g \prod_3^r (g^{2^{j}} \otimes
  Q^j(g^{-1}Qg))  \end{align*}

Using that $Q(g^{-1}Qg) = g^{-2}Q^2(g) + g^{-4}(Qg)^3$,

\bd \psi(x)|_{H_5 \otimes H_{d-5}} = g^2(g^{-1}Qg)^5 \otimes
g^{-1}(Q^2g)^2 \prod_3^r (g^{2^{j}} \otimes
  Q^j(g^{-1}Qg))   \ed

This implies that $5 \in S(x)$, proving the result when $k+1$ is a
power of 2.

\qed
\section{Some Questions}

$\F_p$-cohomology rings of $Rat_3$, $B\beta_{6}$ and $C_3$ are
isomorphic to each other. For $p > 2$, it follows from section 6,
\cite{tot90}. For $p = 2$, we can see this by hand.

$H_*(B\beta_6)$ is the span of monomials of weight 6 in $\F_2[g, Qg,
Q^2g, \cdots]$. And these are $g^6$ (dim 0), $g^4Qg$ (dim 1),
$g^2(Qg)^2$ (dim 2), $g^2Q^2g$, $(Qg)^3$ (both dim 3) and $Q^2gQg$
(dim 4). $H_*(Rat_3)$ is the span of monomials of weight 3 in $\F_2[g,
g^{-1}Qg, Q(g^{-1}Qg), \cdots]$. And these are $g^3$ (dim 0), $gQg$
(dim 1), $g^{-1}(Qg)^2$ (dim 2), $g^{-3}(Qg)^3$, $gQ(g^{-1}Qg)$ (both
dim 3) and $g^{-1}QgQ(g^{-1}Qg)$ (dim 4). It is easy to check from
above that $H_*(B\beta_6)$ is isomorphic as a 
coalgebra to $H_*(Rat_3)$ proving that $Rat_3$, $B\beta_{6}$ and $C_3$
have isomorphic $\F_2$-cohomology rings.

We can also see by hand that the action of the Steenrod algebra on $H_*(Rat_3)$ and $H_*(B\beta_6)$ is the same. Consider $g^2Q^2g \in H_3(B\beta_6)$. Then \begin{align*} Sq_1^*(g^2Q^2g) = g^2(Qg)^2. \end{align*}
The element corresponding to $g^2Q^2g$ in $H_3(Rat_3)$ is $gQ(g^{-1}Qg)$. 
\begin{align*} Sq_1^*(gQ(g^{-1}Qg)) =& Sq_1^*(g^{-1}Q^g) + g^{-3}(Qg)^3 \\ =& g^{-1}(Qg)^2. \end{align*}
$g^2(Qg)^2 \in H_2(B\beta_6)$ corresponds to $g^{-1}(Qg)^2 \in H_2(Rat_3)$. Similarly, by checking for each generator, we can verify that the action of $Sq_1^*$ on $H_*(Rat_3)$ and $H_*(B\beta_6)$ is the same.

It is still unknown if $Rat_3$, $B\beta_{6}$ and $C_3$ have homotopy equivalent
$HF_2$-localizations or not. Also it is still unknown for $p > 2$ and
$k > 2$, if $Rat_k$,
$B\beta_{2k}$ and $C_k$ have homotopy equivalent $H\F_p$-localizations
or not.

Cohen-Shimamoto \cite{cohshi91} have shown that $Rat_2$ and $C_2$ are
not homotopy equivalent to each other by considering natural $\Z$-coverings of these spaces.
\bd \pi_1(Rat_k) \cong \pi_1(C_k)
\cong \Z. \ed

Let $\urat_k$ and $\uc_k$ be the universal covers of $Rat_k$ and $C_k$
respectively. Let $D_2$ be the $\Z/2$-Moore space $S^2 \cup_2 e^3$. It
is known that $C_2$ is stably homotopy equivalent to $S^1 \vee D_2$
\cite{sna74}. Cohen-Shimamoto show that $C_2$ is homotopy equivalent
to $S^1 \vee D_2$. Hence we can precisely calculate \bd H_2(\uc_2; \Z) \cong
\pi_2(\uc_2) \cong \pi_2(C_2) \ed which is infinitely
generated. Whereas $\urat_k$ is homotopy equivalent to $R^{-1}(\{1\})$
for the resultant map $R: Rat_k \ra \C^*$. $R^{-1}(\{1\})$ is a finite
$CW$-complex and hence $H_*(\urat_k; \Z)$ is finitely generated. 

This shows that $\Z$-homologies of $\urat_2$ and $\uc_2$ are not
isomorphic and that \begin{align*} \pi_2(C_2) \ncong \pi_2(Rat_2). \end{align*}. 

Let $\gamma_k$ be the commutator subgroup $[\beta_k, \beta_k]$ of
$\beta_k$. Then there is a short exact sequence \bd 0 \ra \gamma_k \ra
\beta_k \ra  \Z \ra 0. \ed

Let $B\gamma_k$ denote the classifying space of $\gamma_k$. The
$\F_p$-homology of $B\gamma_k$ is calculated in (theorem 4, \cite{cal06}), and it is
finitely generated.

We conjecture that homology of $\uc_k$ is infinitely generated for
many values of $k$ whereas we already know that homologies of $\urat_k$ and $B\gamma_k$ are finitely generated for all $k$.

\def\cprime{$'$} \def\cprime{$'$} \def\cprime{$'$}


\begin{thebibliography}{10}

\bibitem{arn68}
V.~I. Arnol{\cprime}d.
\newblock Braids of algebraic functions and cohomologies of swallowtails.
\newblock {\em Uspehi Mat. Nauk}, 23(4 (142)):247--248, 1968.

\bibitem{bou75}
A.~K. Bousfield.
\newblock The localization of spaces with respect to homology.
\newblock {\em Topology}, 14:133--150, 1975.

\bibitem{boyman88}
C.~P. Boyer and B.~M. Mann.
\newblock Monopoles, nonlinear {$\sigma$} models, and two-fold loop spaces.
\newblock {\em Comm. Math. Phys.}, 115(4):571--594, 1988.

\bibitem{cal06}
F.~Callegaro.
\newblock The homology of the {M}ilnor fiber for classical braid groups.
\newblock {\em Algebr. Geom. Topol.}, 6:1903--1923 (electronic), 2006.

\bibitem{coh78}
F.~R. Cohen.
\newblock Braid orientations and bundles with flat connections.
\newblock {\em Invent. Math.}, 46(2):99--110, 1978.

\bibitem{cohcoh91}
F.~R. Cohen, R.~L. Cohen, B.~M. Mann, and R.~J. Milgram.
\newblock The topology of rational functions and divisors of surfaces.
\newblock {\em Acta Math.}, 166(3-4):163--221, 1991.

\bibitem{cohcoh93}
F.~R. Cohen, R.~L. Cohen, B.~M. Mann, and R.~J. Milgram.
\newblock The homotopy type of rational functions.
\newblock {\em Math. Z.}, 213(1):37--47, 1993.

\bibitem{cohdav88}
F.~R. Cohen, D.~M. Davis, P.~G. Goerss, and M.~E. Mahowald.
\newblock Integral {B}rown-{G}itler spectra.
\newblock {\em Proc. Amer. Math. Soc.}, 103(4):1299--1304, 1988.

\bibitem{cohlad76}
F.~R. Cohen, T.~J. Lada, and J.~P. May.
\newblock {\em The homology of iterated loop spaces}.
\newblock Lecture Notes in Mathematics, Vol. 533. Springer-Verlag, Berlin,
  1976.

\bibitem{cohshi91}
R.~L. Cohen and D.~H. Shimamoto.
\newblock Rational functions, labelled configurations, and {H}ilbert schemes.
\newblock {\em J. London Math. Soc. (2)}, 43(3):509--528, 1991.

\bibitem{fuk70}
D.~B. Fuks.
\newblock Cohomology of the braid group {${\rm mod}\ 2$}.
\newblock {\em Funkcional. Anal. i Prilo\v zen.}, 4(2):62--73, 1970.

\bibitem{hav95}
J.~W. Havlicek.
\newblock The cohomology of holomorphic self-maps of the {R}iemann sphere.
\newblock {\em Math. Z.}, 218(2):179--190, 1995.

\bibitem{seg79}
G.~Segal.
\newblock The topology of spaces of rational functions.
\newblock {\em Acta Math.}, 143(1-2):39--72, 1979.

\bibitem{sna74}
V.~P. Snaith.
\newblock A stable decomposition of {$\Omega \sp{n}S\sp{n}X$}.
\newblock {\em J. London Math. Soc. (2)}, 7:577--583, 1974.

\bibitem{tot90}
B.~Totaro.
\newblock The cohomology ring of the space of rational functions.
\newblock Preprint, MSRI, 1990.

\bibitem{ver99}
V.~V. Vershinin.
\newblock Braid groups and loop spaces.
\newblock {\em Uspekhi Mat. Nauk}, 54(2(326)):3--84, 1999.

\end{thebibliography}
\end{document}